\documentclass[12pt]{article}
\usepackage{amsmath,amssymb,amsthm,amscd,graphics,epsfig}
\usepackage[english]{babel}

\newtheorem{theorem}{Theorem}
\newtheorem{lemma}[theorem]{Lemma}

\newtheorem{corollary}[theorem]{Corollary}
\newtheorem{definition}[theorem]{Definition}

\begin{document}

\title{Poncelet pairs and the  Twist Map associated to the Poncelet Billiard}
\author{ A. O. Lopes\footnote{ Instituto de Matem\'atica, UFRGS, Porto
Alegre 91501-970, Brasil. Partially supported by CNPq, Instituto do
Milenio, PRONEX - Sistemas Din\^amicos, benefici\'ario de
aux\'{\i}lio financeiro CAPES - PROCAD.}\and M. Sebastiani
\footnote{ Instituto de Matem\'atica, UFRGS, Porto Alegre 91501-970,
Brasil.} } \maketitle

\begin{abstract}

Consider a fixed differentiable
curve  $K$ (which is the boundary of a convex domain) and a family indexed by $\lambda\in [0,1]$ of variables differentiable
curves $L$ such that are the boundary of convex domains $A_\lambda$. Suppose that for $\lambda_1<\lambda_2$ we have
$A_{\lambda_1}\subset A_{\lambda_2}$. Then, the number of $n$-Poncelet pairs is given by  $\frac{e
(n)}{2}$, where $e(n)$ is the number of natural numbers $m$
smaller than $n$ and which satisfies mcd $\,(m,n)=1$. The curve
$K$ do not have to be part of the family.

In order to show this result we consider an associated billiard
transformation and a twist map which preserves area.

We use Aubry-Mather theory  and the rotation number of invariant
curves to obtain our main result.

In the last section we estimate the derivative of the rotation number of a general twist map
using some properties of the continued fraction expansion .

\end{abstract}

\newpage

\section{Introduction - The Poncelet Billiard}

All results here are for a fixed differentiable
curve  $K$ (which is the boundary of a convex domain) and a family indexed by $\lambda\in [0,1]$ of variables differentiable $C^{\infty}$ 
curves $L$ such that are the boundary of convex domains $A_\lambda$. Suppose that for $\lambda_1<\lambda_2$ we have
$A_{\lambda_1}\subset A_{\lambda_2}$.

Some explicit results (not for the Poncelet counting problem) for billiards in the case the curves are ellipses were obtained in \cite{c1}. Our reasoning apply for a more general family of curves. We refer the reader to \cite{c2}
for some applications of the Poncelet billiard problem.

In order to simplify the exposition we will refer to circles instead of curves. Some computations are explicit for circles but our reasoning applies to the setting we just described above. The estimates of Proposition 2 bellow are 
not contained in previous results which analyze the Poncelet 
problem.

Let $K$ a circle as in figure 1, and also $L$ another circle
interior to $K$ with $A, A'$ variable points in $K$.

Consider the Poncelet transformation associated to such pair of
circles $K$ and $L$ and the corresponding image of the point $A$
being $A'$, and then the image of $A'$ is $A''$ and so on... (see
figure 1).  A nice description of the problem appears in \cite{K}and \cite{T}.

$T, C, o$ are shown in figure 2. Here $R$  is the radius the
radius of the circle $K$ and $C$ is the center of $K$.

We will  use the variables $\theta$ and $\varphi$ to describe the
point $A$ and his future hit as in figure 1 and 2. The point $A'$
has coordinates $\theta'$ and $ \varphi'$.

The variable points $B,B'$ on the $x$-axis are also described in
figure 2.

Denote $ \theta'=G_1(\theta, \varphi)$ and $ \varphi'=G_2(\theta,
\varphi).$ The transformation $G=(G_1,G_2)$ can be consider also
as a transformation of $(A, r)\to (A', r')$, where $r$ and $r'$
are tangent lines to the circles as in figure 1.

Note that once $A$ is fixed, if we take a variable  line $r$, this
means we are considering different circles $L$. The circle $K$ is
considered fixed in our setting.

The reader can find several references about the Poncelet Billiard
in [1], [3], [4], [6], [7], [8], [10].

We say $K, L$ is a $n$-Poncelet pair if there exists an $A$ such
that the successive iteration of the procedure described above,
after $n$ steps, returns to $A$.

The existence of a $n$-periodic point $(A,r)$ for $G$ (for a
certain pair of circles $K,L$) is equivalent to the existence of a
$n$-Poncelet pair. We would like to count the number of possible
Poncelet pairs of order $n$.

It is easy to obtain the analytic expression
$$G(\theta,\varphi)=( 2 \, \varphi - \theta + \pi,\, 3\, \varphi- 2 \theta  - B (\theta') +
\pi), $$ where $B(\theta') = 2\, \arctan   (  \frac{C \, \sin
(\theta ')  }{R + C \cos(\theta')}).$

\begin{figure}
\begin{center}
\epsfig{file=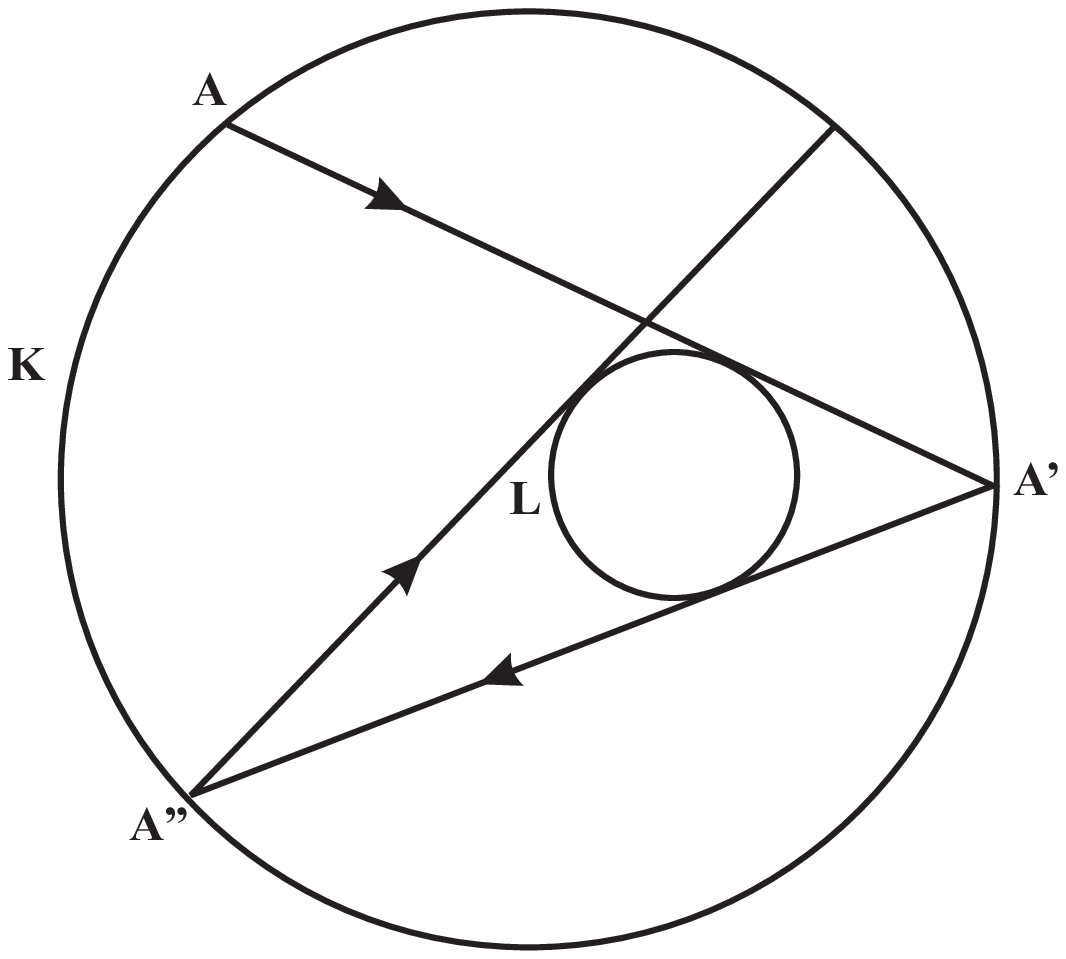, width=8cm} \caption{ }
\end{center}
\end{figure}

We know that $r$ is a periodic function of $\varphi$ of period
$\pi$. We change coordinates to get $x= \frac{\theta}{ 2 \, \pi}$
and $y= \frac{\varphi}{  \, \pi}$.

Therefore, denoting by $f$ the transformation $G$ in the new
coordinates
$$ (x',y')=f(x,y) = (\,y\,-x\,+\,1/2, \, 3 \, y - 4 x + Z( y-x +
1/2\,)\, +\,1\, ),$$

where $Z(s) = -\, B( \,2 \, \pi \,s).$ The function $Z$ is
periodic of period $1$.

Therefore $f=(f_1,f_2)$ is a twist map on the torus (see
[2],[5],[9] for references), because

$$ 1) \,\,\,\,\,\frac{\partial f_1(x,y) }{\partial y}=\frac{\partial x '}{\partial y} =1>0,$$

and, also preserves area

$$ 2) \,\,\,\,\, \mbox{Jac}\, D\, f \,
=\,\mbox{Det}\, \begin{pmatrix} -1\,\, &  1\\ -4 - Z'(y-x+1/2)\,\,
&  3 + Z'(y-x+1/2)
\end{pmatrix} =1.$$

One can write $y = x'+ x - 1/2$, and therefore the generating
potential $h(x,x')$ for such $f$ is
$$ h(x,x') \,=\, -\, x \,x' -\, \frac{x^2- x}{2}+ \frac{3\, (x')^2 -
x'}{2}+ H(x'),$$ where $H$ is such that $H'(x') = Z(x').$

The twist map is not exact.

After this brief introduction we will present in next section a
sequence of results that will be used in the last section to show
our main theorem:

\begin{theorem}

For a fixed circle $K$ and for a family of variables concentric
circles $L$, the number of $n$-Poncelet pairs is $\frac{e
(n)}{2}$, where $e(n)$ is the number of natural numbers $m$
smaller than $n$ and which satisfies mcd $\,(m,n)=1$.

\end{theorem}

The general problem for conics is the following: given a conic $K$
and a pencil of conics ${\cal L}$, obtain the number of conics $L
\in {\cal L}$ such that $K,L$ is a $n$-Poncelet pair. This problem
was solved in [11] for the complex plane for a generic pencil  but
under the condition that $K \in {\cal L}$. In our analysis we have
that $K$ is not in ${\cal L}$.
\vspace{0.3cm}

Consider a family $g_t : \mathbb{R} \to \mathbb{R}$, where $a\leq
t\leq b$, of monotonous increasing homeomorphisms, such that
$$ g_t (x+1) = g_t (x) \,+\,1,\, x\in \mathbb{R}, \,a\leq
t\leq b.$$

From the identification $S^1 = \frac{\mathbb{R}}{\mathbb{Z}},$
from the family $g_t$, we obtain another one $f_t: S^1 \to S^1$,
where $a\leq t\leq b$, of homeomorphisms of the circle which
preserve orientation.

We know that is well defined
$$ r(t)\, = \, \lim_{k \to \infty} \frac{g_t^k (x) - x}{k},\,\,
a\leq t\leq b,$$ and the value $r(t)$  indeed is independent of
$x$.

We are interested in properties of the function $ r(t)$ and of the
rotation number function $\rho(t) = r(t) $ (mod $1) \in
\frac{\mathbb{R}}{\mathbb{Z}}$.

In the last section we show

{\bf Proposition 1:} The function $r:[a,b] \to \mathbb{R}$ is
continuous monotonous increasing. Moreover, if $   a\leq t_1 < t_2
\leq b$ and $r(t_1)\notin \mathbb{Q} $ or $r(t_2)\notin
\mathbb{Q}$, then $r(t_1) < r(t_2)$.

And also,

{\bf Proposition 2:} There exist a set of full measure $X \subset
\mathbb{R}$ (that is, $ \mathbb{R}- X$ has Lebesgue measure zero)
and such that if $\tau \in (a,b)$ and $r(\tau) \in X$, then
$$ \limsup_{t_2 \to \tau^- , \, t_1 \to \tau^+}\, \frac{r(t_2)-r(t_1) }{  (t_2\,-\,t_1)^2
}\geq \frac{m^2}{e^{2 \, F} \, (1 + e^{2\, F} )^2}.$$

\begin{figure}
\begin{center}
\epsfig{file= 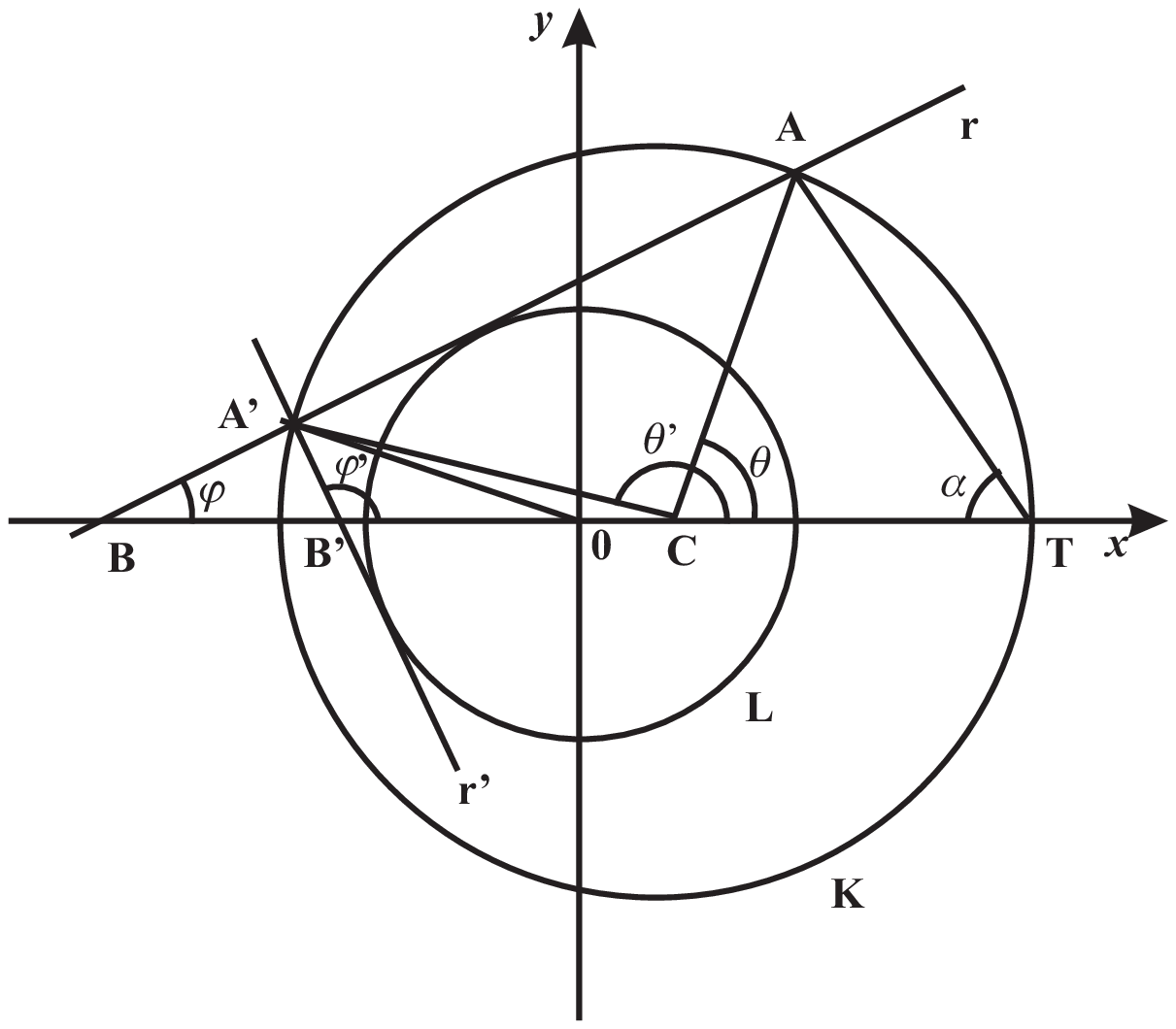  , width=8cm} \caption{ }
\end{center}
\end{figure}

\section{Results for  Twist Maps }

Consider $f: S^1 \to S^1$, where $S^1=
\frac{\mathbb{R}}{\mathbb{Z}}$, a circle homeomorphism which
preserves orientation. Denote by $g:\mathbb{R} \to \mathbb{R}$ a
lifting of $f$.

Denote by $r(g)$ the limit
$$ \lim_{k \to \infty} \frac{g^k (x) - x}{k}= r(g),$$
which is independent of $x\in \mathbb{R}$ (see [9] chapter 11).

The rotation number of $f$, denoted by $\rho(f)$, is the number
$r(g)$ (mod $1) \in S^1$.

General references for the rotation number and Twist Maps are [2],
[5], [9].

\begin{lemma} For $f$ and its lifting $g$:

a) $\rho (f)$  is an invariant of conjugation (by an homeomorphism
which preserves orientation)  for $f$,

b) if $\rho(f)= \frac{p}{q}$ (mod $1$), where $q>0$ and $p,q$ are
relatively prime, then there exists a $x_0$ such that $g^q(x_0)
=x_0 +p$ (see [9]), and the orbit of $x_0$ (mod 1) is periodic of
period exactly $q$.
\end{lemma}

The following Lemma is also well known:

\begin{lemma}

Consider $g,h: \mathbb{R} \to \mathbb{R}$ continuous and strictly
increasing such that $g(x+1)= g(x) + 1$ and  $h(x+1)= h(x) +1$,
$\forall x \in \mathbb{R}$.

Then,

a) If $h(x) \leq g(x),$ $\forall x \in \mathbb{R}$, then $r(h)\leq
r(g)$.

b) $g \to r(g)$ is continuous in the $C^0$ topology.

c) If $h(x) < g(x),$ $\forall x \in \mathbb{R}$, and if $r(h)$ or
$r(g)$ is irrational,  then $r(h)<  r(g)$.

\end{lemma}

Now, consider the cylinder $X= S^1 \times \mathbb{R} =
\frac{\mathbb{R}}{\mathbb{Z}} \times \mathbb{R}$ and the natural
projection $\pi: \mathbb{R}^2 \to X$.

\begin{definition}

The vector field $u$ over  the cylinder $X$ is called the
ascendent vector field, if $u$ is such that  $u(a , b)=$ the
tangent vector to the curve $t \to (a , t\,+\,  b)$, in the point
$t=0$.

\end{definition}

We denote by $ \pi_1: X \to S^1$ and $\pi_2: X \to \mathbb{R}$ the
natural projections.

\vspace{0.3cm}

We will consider here diffeomorphisms of class $C^1$, $f:X \to X$,
such that

a) $f$ preserves orientation and the ends of $X$,

b) $f$ preserves area of a Riemannian metric, and

c) for each point $x\in X$ the vectors $u$ and $df_x (u)$ are
independent and coherent with the natural orientation.
\vspace{0.3cm}

Consider $F:\mathbb{R}^2 \to \mathbb{R}^2$, $F(x,y)=(F_1(x,y),
F_2(x,y))$, a lifting of such $f$.

Then $F$ satisfies the conditions

1) $F_1(x+1,y)= F_1(x,y) + 1$,

2) $F_2(x+1,y) = F_2 (x,y)$,

3) for each $x$ , we have $\lim_{y \to \infty} F_2 (x,y)= \infty$,

4) $F$ preserves area of a Riemannian metric which is invariant by
translation $(x,y)\to (x+1,y)$,

5) $\frac{\partial F_1}{\partial y}\,(x,y)\,>\,0$ (the twist
condition). \vspace{0.3cm}

We point out that by condition 5) we have that given $(x_0,y_0)$ and
$(x_0, y_1)$ in $\mathbb{R}^2$, such that, $y_1>y_0$, then
$F_1(x_0,y_1) > F_1 ( x_0, y_0).$

\vspace{0.3cm}

\begin{definition}

Consider $\phi: \mathbb{R} \to \mathbb{R}$, a continuous function of
period $1$, and the curve $\Gamma\subset X$ with parametrization
$x\to \pi (x, \phi(x))$, $0\leq x\leq1.$
 We say $\Gamma$ is an invariant rotational circle for $f$ if
 $f(\Gamma)=\Gamma.$

\end{definition}

We point out that $ \Gamma$ is an oriented circle and $\pi^{-1}
(\Gamma)$ is the graph of $\phi$. Moreover $f:\Gamma\to \Gamma$
preserves the orientation of $\Gamma$. We denote by $\rho(\Gamma)$
the rotation number of $\rho (f|\Gamma).$ The transformation $F$
also preserves the graph of $\phi$ and its orientation.

The usual definition of invariant rotational circle is more
general ([2] chapter 3, definition 11) but a theorem due to
Birkhoff ([2], 3.1) shows that, under quite general conditions,
this definition coincides with our definition 5. The rotational
invariant curve we are going to consider here (a subset of the
torus) is a set of positions $q$ on the circle $K$ and respective
tangent line tangent to the small circle $L$ (passing through
$q,p$) as shown in figure 3.

Consider a $g$ associated to $\Gamma$ by $g(x) = F_1(x, \phi(x)),
\forall x \in \mathbb{R}$. The function $g$ is continuous strictly
increasing and satisfies
$$ g(x+1) = F_1(x+1, \phi(x+1))= F_1 (x+1, \phi(x)) = g(x) + 1$$

We claim that
$$\rho(\Gamma) = r(g) \, (\text{mod} 1).$$

Indeed, by definition $\pi_1| \Gamma : \Gamma \to S^1$ is a
homeomorphism which preserves orientation. Consider $\varphi= (
\pi_1 | \Gamma) \circ (f | \Gamma) \circ (\pi_1 | \Gamma)^{-1}.$
 Then $\varphi$ preserves orientation and $\rho (\Gamma) = \rho
 (\varphi)$.

By the other side
$$ \pi_1 \circ \pi (x,0)= x \, (\text{mod} 1), \forall x \in \mathbb{R}.$$

As
$$\varphi(\pi_1 (\pi(x,0)))= \varphi (\pi_1 (\pi(x, \phi(x))))=$$
$$ \pi_1(f(\pi(x, \phi(x))))= \pi_1(\pi(F(x,\phi(x))))= \pi_1\circ
\pi (g(x),0),$$ then $g$ is a lifting of $\varphi$. Then, $\rho
(\varphi)=r(g)\, ($mod $1)$ and  the claim is true.

\begin{theorem}

Consider $\Gamma_1$ and $\Gamma_2$ two rotational invariant circles
associated respectively to $\phi_1$ and $\phi_2$.

Suppose $\phi_1(x)< \phi_2 (x), \forall x \in \mathbb{R}$. Denote
$$ g_1(x)=F_1 (x, \phi_1 (x)), \,\,\,     g_2(x)=  F_1 (x, \phi_2 (x)), \,\,\forall x \in \mathbb{R}.  $$

Then
$$ r(g_1) < r(g_2).$$

\end{theorem}

{\bf Proof:} Following [2] 3.3,  from $\phi_1<\phi_2$ we have that
that $g_1(x) < g_2 (x), \forall x \in \mathbb{R}.$ Therefore,
$r(g_1) \leq r( g_2)$. If $r(g_1)$ or $r(g_2)$ is irrational, then
the claim is true (Lemma 3).

Suppose, that $r(g_1)=r(g_2)= \frac{p}{q},$ where $q>0$ and  $p,q$
are relatively prime.

Consider the translation $T_n : \mathbb{R}^2 \to \mathbb{R}^2$
given by $ T_n(x,y) = (x+n,y)$ and the diffeomorphism $G=T_{-p}
\circ F^q : \mathbb{R}^2 \to \mathbb{R}^2$. We will show that
$G(A)$, for a certain region $A$, is strictly inside the set  $A$
and we will get a contradiction with the fact that $G$ preserves
area.

We denote $ m'= G(m)$, $m \in \mathbb{R}^2$.

Note that $G$ preserves the oriented graphs of $\phi_1$ and
$\phi_2$.

By Lemma 2 b) there exists $x_1,x_2\in \mathbb{R}$ such that
$$ g_1^q (x_1)= x_1 + p, \, \, \, g_2^q (x_2) = x_2 + p.$$

Without lost of generality suppose $ x_2>x_1$.

We denote
$$ M_1=( x_1, \phi_1 (x_1)),  M_2=( x_2, \phi_2 (x_2)),
M_3=( x_1, \phi_2 (x_1)), M_4=( x_2, \phi_1 (x_1)).$$

As
$$ F^q (x, \phi_1 (x) )= (g^q (x), \phi_1 (g^q(x))= (x+ p, \phi_1
(x+p)) =(x+p, \phi_1 (x)),$$ then $ M_1' = G(M_1)= G(x_1, \phi_1
(x)) =(x_1, \phi_1(x))=M_1$.

In the same way $M_2'= M_2$

By the twist condition, $M_3 '$ follows $M_3$ in the orientation
of the graph of $ \phi_2$ (see figure 4). As $M_2$ is fixed, then
$M_3'$ precedes in order to $M_2$. That is, $M_3'$ is interior to
the arc (which is part of the graph of $\phi_2$) connecting $M_3$
and $M_2$ (in this order as shown in figure 4). In the same way,
$M_4'$ is interior to the arc (which is part of the graph of
$\phi_1$) connecting $M_1$ and $M_4$ (in this order).

Furthermore, by the twist condition, the image of the linear
interval $L_1= [M_1,M_3]$ by $G$ is situated to the right of  the
line defined by $M_1$ and $M_3$. By the same reason, the image of
the linear interval $L_2= [M_4,M_2]$ by $G$ is situated to the left
of the line defined by $M_4$ and $M_2$.

Finally, as the graphs of $\phi_1$ and $\phi_2$ are invariant by
$G$, then it follows from above that the  region $A$ delimited by
them (and the lines $M_1,M_3$ and $M_2,M_4$) is mapped strictly
inside $A$ by the transformation $G$.

This is a contradiction with the fact that $G$ preserves area.

\qed

\begin{corollary}

Let $\phi_t: \mathbb{R} \to \mathbb{R} $, $t\in [a,b]$, be a family
of continuous periodic functions of period $1$. Suppose that
$\phi_t(x)$ is continuous as a function of $(t,x)$ and that, for
each $t$, the transformation $\phi_t$ defines a rotational invariant
circle $\Gamma_t$, such that they are disjoint two by two.

Then there exists a continuous strictly monotone function $r: [a,b]
\to \mathbb{R}$ such that
$$ \rho(\Gamma_t) = r(t) \, (\text{mod }\, 1), \,\,\,\,\forall t \in [a,b].$$

\end{corollary}

{\bf Proof:} By hypothesis, for each $x \in \mathbb{R}$, the
function $t \to \phi_t (x)$ is injective. Therefore, it is either
strictly increasing or strictly decreasing. If it is increasing for
some point $x$, it  can not be decreasing for other points.

\begin{figure}
\begin{center}
\epsfig{file=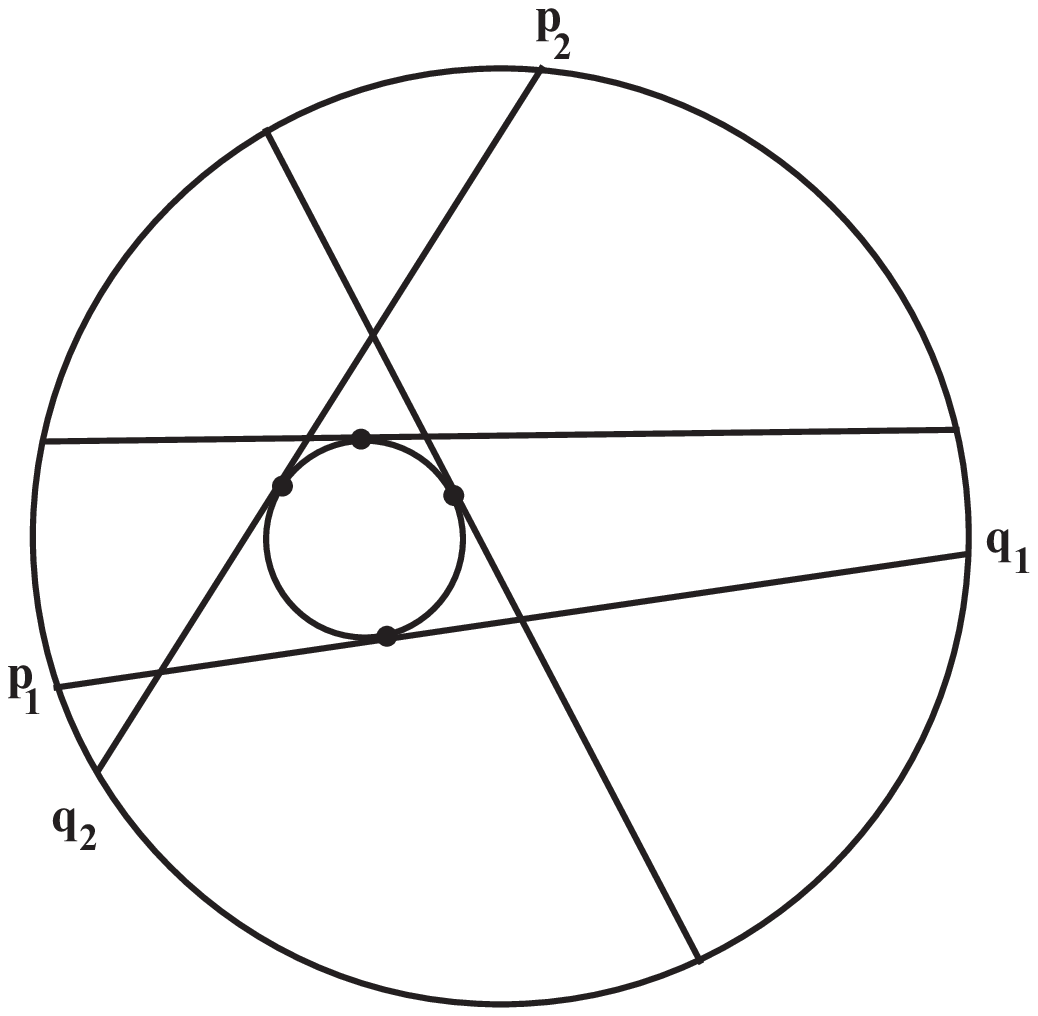, width=8cm} \caption{ }
\end{center}
\end{figure}

Without lost of generality, suppose that $t \to \phi_t (x)$ is
strictly increasing for all $x\in \mathbb{R}$.

Denote $g_t(x)= F_1 (x, \phi_t (x))$. Then, by the claim after
Definition 5, for all $t\in [a,b]$,
$$ \rho (\Gamma_t) = r(g_t) \, (\text{mod} 1).$$

Therefore, take $r(t)=r(g_t)$. Then, $r$ is continuous (see [2]
[9]) and strictly increasing (by Theorem 6).

\qed

Consider now the torus $Y= S^1 \times S^1$ with the natural
orientation and let $\pi_1, \pi_2: Y \to S^1$ be  the canonical
projections. Denote by $v$ the unitary vector field  which is
tangent to the fibers of $\pi_1$ (compatible with the orientation of
$S^1$).

Consider a $C^1$  diffeomorphism $f:Y\to Y$ such that

a) $f$ preserves orientation of $Y$,

b) $f$ preserves area of a Riemannian metric on $Y$,

c) $v, df(v)$ are independent in each point and they are compatible
(in this order) with the orientation of $Y$ (twist condition).

\begin{definition}

Consider $\phi: S^1 \to S^1$, a continuous function and  the
oriented curve $\Gamma\subset X$, with parametrization $x\to \pi
(x, \phi(x))$, $x\in S^1.$
 We say $\Gamma$ is an invariant rotational circle for $f$, if
 $f(\Gamma)=\Gamma$ and $f|\Gamma$ preserves orientation.

\end{definition}

We point out that $\rho (f|\Gamma)$ is well defined and we denote
$\rho(\Gamma)= \rho ( f | \Gamma).$

\begin{corollary}

Let $\phi_t: S^1 \to S^1 $, $t\in [a,b]$, be a family of continuous
functions. Suppose that $\phi_t(x)$ is continuous as a function of
$(t,x)$ and that for each $t$ the transformation $\phi_t$ defines a
rotational invariant circle $\Gamma_t$, such that they are disjoint
two by two.

Then there exists a continuous strictly monotone function $r: [a,b]
\to \mathbb{R}$ such that
$$ \rho(\Gamma_t) = r(t) \, (\text{mod } 1), \,\,\,\,\forall t \in [a,b].$$

\end{corollary}

{\bf Proof:}  Denote $Y_0 = Y- \Gamma_b$, then there exists a
diffeomorphism $\theta : X \to Y_0$ which preserves orientation,
compatible with the projections $\pi_1 :X \to S^1$, $\pi_1 | Y_0:
Y_0 \to S^1$, and such that preserves the orientation  of the fibers
of these projections.

The transformation $f_0= \theta^{-1} \circ f \circ \theta : X \to
X$ is a diffeomorphism of the cylinder $X$ and which satisfies the
properties of $f:X \to X$ described just after Definition 4 ($f_0$
preserves the ends of $X$ because it preserves the orientation of
$Y$ and the orientation of $\Gamma_b$).

For each $t\in[a,b],$ we have that $\theta^{-1} (\Gamma_t)$ is a
rotational invariant circle of $f_0$. In the same way as in
Corollary 7, there exists $r[a,b) \to \mathbb{R}$, which is
continuous, strictly monotone and such that
$$ \rho(\Gamma_t) =r(t)\, (\text{mod} \, 1),\,\,\, t \in [a,b).$$

Note that $\pi_1 | \Gamma_t: \Gamma_t \to S^1$ is a diffeomorphism
which preserves orientation.

By conjugacy,
$$ \rho( \Gamma_t) = \rho( (\pi_1| \Gamma_t) \circ (f| \Gamma_t)
\circ (\pi_1| \Gamma_t)^{-1}\,).$$

Then, by Lemma 3 b), we have that
$$ \lim_{t \to b^{-}} \rho (\Gamma_t) = \rho (\Gamma_b).$$

From this follows that
$$ \lim_{t \to b^{-}} r (t)=l,$$
is finite. We define $r(b)=l$ and
$$ \rho (\Gamma_b) = r(b) \, (\text{ mod }1).$$

\qed

\section{The main  result}
\vspace{0.3cm}

Using the previous notation for the Poncelet billiard, denote by $Y$
the torus
$$ Y= \{ (A, r): A \in K,\, r \,\text{ a line by}\, A\}= K\times S^1$$

The transformation $f:Y \to Y$ associated to the Poncelet Billiard
was defined before in section 1. $f$ satisfies the hypothesis
described in section 2.

\begin{figure}
\begin{center}
\epsfig{file=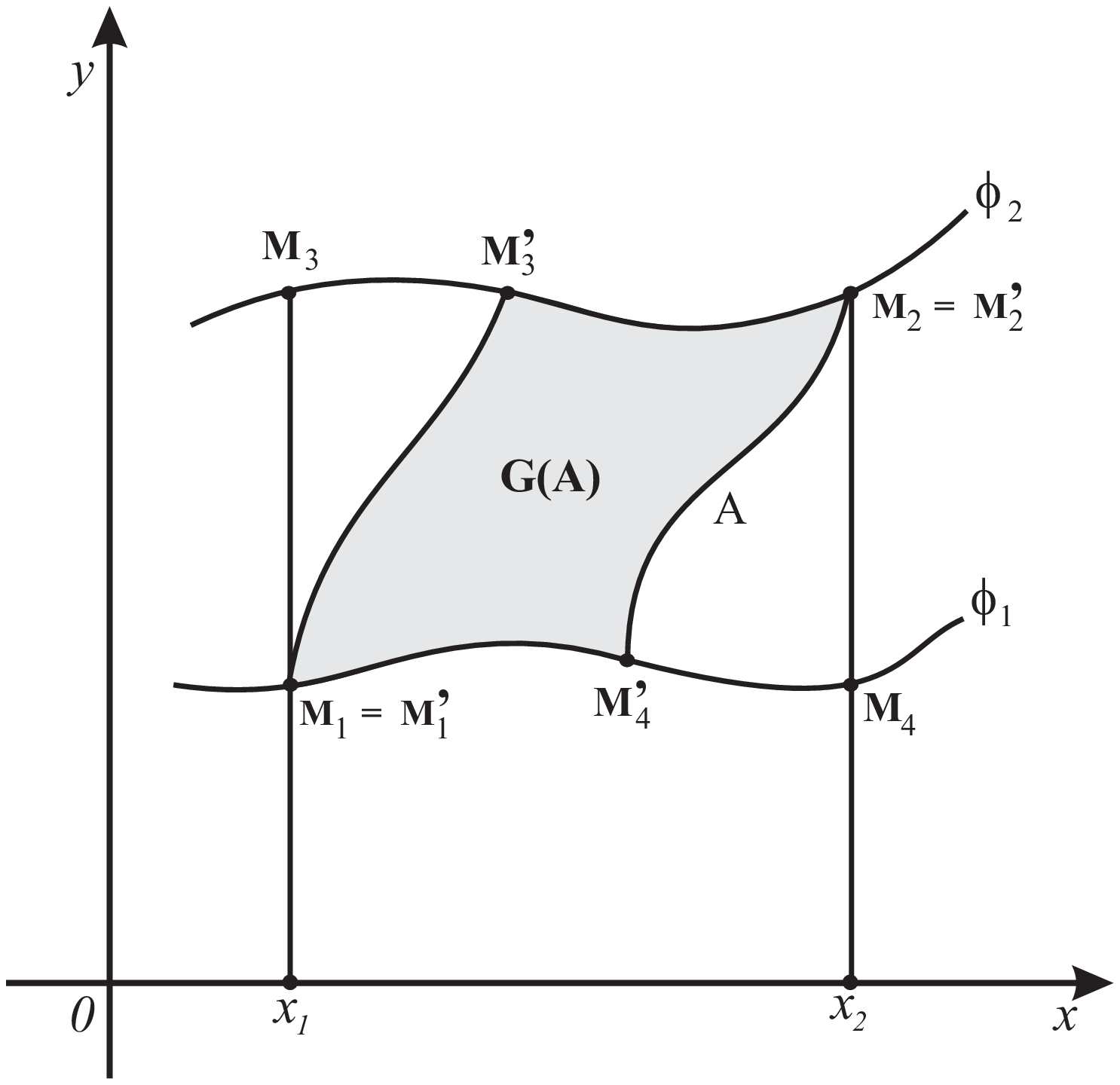, width=8cm} \caption{ }
\end{center}
\end{figure}

Consider a fixed circle $L$ of center $0$ and radius $t$, where $t
\in [0, R-c]$.  To each point $A\in K=S^1$, one can associate a
tangent $r$ to $L$ passing by $A$, and which has the property that
$L$ is on the left of $r$ when oriented from $A$ to $A'$.

In figure 3 we show the set of tangents to $L $ defining the
corresponding rotational invariant circle $\Gamma_t$.

For a fixed circle $L$, we define in such way a transformation
$\phi_t : S^1 \to S^1$, which depends continuously on $t$. These
$\phi_t$ define a family of rotational invariant circles $\Gamma_t$,
which are, two by two, disjoint.

By corollary 9, there exists $r:[0, R-c] \to \mathbb{R},$
continuous, strictly monotonous and such that
$$ \rho(\Gamma_t) = r(t) \, (\text{mod} 1).$$

As $f| \Gamma_{R-c} $ has a fixed point, then $\rho (\Gamma_{R-c})
\,=\,0\, ($ mod $1).$

As a continuous limit of the figure 1 when $t\to 0$, one can see
that $f| \Gamma_0$ has period $2$, and therefore
$\rho(\Gamma_0)\,=\,1/2\, ($ mod $1)$.

Therefore, $r(0)=\frac{1}{2} + n, $ and $r(R-c)=m$,  where $n,m \in
\mathbb{Z}$.

We can suppose that $r(t)$ is such that $ r(0)=1/2$, and $r(R-c)=k
\in \mathbb{Z}$.

We point out that $f| \Gamma_t$ has no fixed points, if $t\neq
(R-c)$. Then, $r(t)$ is not in $\mathbb{Z}$,  if, $t\neq (R-c)$
(Lemma 2 b). The conclusion is that $k=0$ or $k=1$. Therefore, $r$
is continuous, strictly increasing  with values in $[1/2,1]$ or
strictly decreasing  with values in $[0,1/2]$.

For each $q=3,4,5,...$, denote $e(q)$ the number of natural numbers
$n<q$, which are relatively prime with $q$. Therefore, in any case,
for a fixed $q$, there exactly $\frac{e(q)}{2}$  numbers $n$ in such
way $\frac{n}{q}$  is attained in the image of $r$.

From this and Lemma 2 b), it follows:

\begin{lemma} For each natural number $q=1,2,3,4,..$, there are
$\frac{e(q)}{2}$ rotational invariant circles $\Gamma_t$, such that
$f|\Gamma_t$ has a periodic orbit of period exactly $q$.

\end{lemma}

By Poncelet theorem [6] and [7], $f|\Gamma_t$ has a periodic orbit
of period $q$, if and only if, $f|\Gamma_t$ is a $q$-periodic
transformation. In this case the corresponding pair $K,L$ is
called a Poncelet pair.

From this follows  theorem 1.

\vspace{0.3cm}

\section{A second order estimate of the derivative of twist maps}

The results of this sections are for a more general class of twist
maps.

Consider a family $g_t : \mathbb{R} \to \mathbb{R}$, where $a\leq
t\leq b$, of monotonous increasing homeomorphisms, such that
$$ g_t (x+1) = g_t (x) \,+\,1,\, x\in \mathbb{R}, \,a\leq
t\leq b.$$

From the identification $S^1 = \frac{\mathbb{R}}{\mathbb{Z}},$
from the family $g_t$, we obtain another one $f_t: S^1 \to S^1$,
where $a\leq t\leq b$, of homeomorphisms of the circle which
preserve orientation.

We know that is well defined
$$ r(t)\, = \, \lim_{k \to \infty} \frac{g_t^k (x) - x}{k},\,\,
a\leq t\leq b,$$ and the value $r(t)$  indeed is independent of
$x$.

We are interested in properties of the function $ r(t)$ and of the
rotation number function $\rho(t) = r(t) $ (mod $1) \in
\frac{\mathbb{R}}{\mathbb{Z}}$.

We assume just the {\bf twist condition}: there exists
$\frac{\partial g_t (x)}{\partial t}$ and is continuous and
positive for all $x \in \mathbb{R}$, $t \in [a,b]$.

As $\frac{\partial g_t (x)}{\partial t} $ is periodic in $x$, we
get
$$ m= \inf_{ a\leq t\leq b,\, x \in \mathbb{R}   }  \,   \, \frac{\partial g_t (x)}{\partial t}      \,> \,0.$$

For a given irrational value $x$ consider its development in
continuous fraction and we call approximations of $x$ the
successive truncations of this infinite expansion.

We denote by defect approximation one which is smaller than $x$
and by excess approximation one which is larger than $x$.

We refer the reader to \cite{DK} for general properties of
continuous fraction expansion.

We denote by $F$ the sum of the inverses of the Fibonacci
sequence:
$$ F = \frac{1}{1}  + \frac{1}{1}  +\frac{1}{2}  +\frac{1}{3}  +\frac{1}{5}
+  \frac{1}{8} + ...$$

Consider two monotonous increasing homeomorphisms $g_1,g_2:
\mathbb{R} \to \mathbb{R}$  such that
$$ g_1 ( x+1) = g_1 (x) +1
\, ,\, g_2 ( x+1) = g_2 (x) +1 \, ,\, x \in \mathbb{R}.$$

Denote
$$ r_1 = \lim_{k \to \infty}  \frac{g_1^k (x) - x}{k},\,r_2 = \lim_{k \to \infty}  \frac{g_2^k (x) - x}{k}
.$$

Suppose $g_1(x) < g_2(x),\, \forall x \in \mathbb{R}$. As $g_1(x)
- g_2 ( x) $ is periodic and continuous, we have that
$$ \alpha = \inf_{x \in \mathbb{R} } \, (g_2(x) -g_1 (x))\,
>\,0.$$

{\bf Lemma 1.}

a) $r_1 \leq r_2$,

b) Suppose that $r_1$ is irrational and denote by $\frac{p}{q}$ an
excess continued fraction approximation of $r_1$ (that is,
$\frac{p}{q}\geq r_1$) such that $q> \frac{1}{\alpha}.$ Then,
$$
r_1\,< \frac{p}{q} \, \leq r_2.$$

c) suppose $r_2$ is irrational and denote by $\frac{p'}{q'}$ an
excess continued fraction approximation of $r_1$ (that is,
$\frac{p'}{q'}\geq r_1$) such that $q'> \frac{1}{\alpha}.$ Then,
$$
r_1\,\leq \frac{p'}{q'} \, < r_2.$$

For a proof of the lemma see section 11.1 in \cite{KH}

We need some results about approximation by truncation of
continued fractions.

{\bf Lemma 2.} For each irrational $x$ consider $
\frac{p_n(x)}{q_n (x)}$ , $ n=1,2,3,...$, the $n$-truncation of
the continuous fractional expansion of $x$. Given $0< \epsilon<1$,
then for Lebesgue almost every irrational $x\in \mathbb{R}$, there
exists a sequence $n_1<n_2<n_3,..$ of natural numbers (which
depend on $x$ and $\epsilon$)  such that
$$ \frac{\,2 \, e^{-2 \, F} }{1+ \epsilon } < \frac{q_{n_k +1} (x)       }{q_{n_k}(x)     }  <
\frac{\,2 \, e^{2 \, F}   }{ 1-\epsilon}, \, \,\,k =1,2,3,...$$

{\bf Proof:} We can suppose $x\in (0,1)$. Denote $T :[0,1)\to
[0,1),$  the Gauss transformation $T(x) = \frac{1}{x}- [
\frac{1}{x} ],$ if $x\neq 0,$ and $T(0)=0$.

For each irrational number $x$, from section 3.5 in \cite{DK}
$$ - \log q_n (x)\,=\, \log(x) + \log (T(x)) + ...+ \log (T^{n-1}
(x)) + R(n,x), \,\, n=1,2,3,...,$$ where $|R(n,x)|\leq F$, for all
$x$ and $n=1,2,3,...$

Therefore,
$$ \log ( \frac{q_{n+1} (x)}{q_n (x)}) = -\log T^n (x) + R(n,x) - R(n+1,x),$$
and finally
$$      - 2\, F - \log T^n (x)    \,\leq\,               \log ( \frac{q_{n+1} (x)}{q_n
(x)})\,\leq\, 2\, F - \log T^n (x) .$$

As $T$ is ergodic \cite{DK}, for Lebesgue almost every irrational
$x\in (0,1)$, we have that $\frac{1}{2},$ is a limit of a certain
subsequence of $T^n (x)$. Therefore, there exists natural numbers
$n_1<n_2<n_3<...$, such that
$$\frac{1}{2} - \frac{\epsilon }{2}< T^{n_k} (x) <\frac{1}{2} + \frac{\epsilon
}{2},\,\, k=1,2,3,....$$

From this we get
$$ - 2\, F - \log (\frac{1}{2} + \frac{\epsilon }{2}  ) \,<\, \log ( \frac{q_{n_k+1}
(x)}{q_{n_k} (x)})\,<\,   2\, F - \log (\frac{1}{2} -
\frac{\epsilon }{2}  ) .$$

\qed

{\bf Corollary 1:}  Given $0<\epsilon<1,$ for almost irrational
$x$ there exists an infinite number of excess approximations
$\frac{p}{q},$ and infinite number of excessive approximations
$\frac{p'}{q'},$ such that
$$ \frac{\,2 \, e^{-2 \, F} }{1+ \epsilon} < \frac{q'}{q}   <
\frac{\,2 \, e^{2 \, F}   }{1-\epsilon}, \, \,\mbox{or}\,\,
\frac{\,2 \, e^{-2 \, F} }{1+ \epsilon} < \frac{q}{q'}   <
\frac{\,2 \, e^{2 \, F}   }{1-\epsilon}.
$$

{\bf Proof:}  In the Lemma 2 above take
$$p=p_{n_k}(x),\,
q=q_{n_k}(x), \,p'=p_{n_k+1}(x), \,q'=q_{n_k+1}(x),$$ if $n_k$ is
odd and the opposite if $n_k$ is even.

\qed

{\bf Corollary 2:}  Given $0<\epsilon<1,$ denote
$$ K_\epsilon= \frac{1- \epsilon}{e^{2\, F} \, (1 + (1+ \epsilon) \, e^{2 \,F})^2}.$$

Then, for almost every irrational $x\in \mathbb{R}$ there exists
infinite excess approximations $\frac{p}{q}$ and infinite defect
approximations $\frac{p'}{q'}$ of x such that:

$$  \frac{p}{q}   \,-\,   \frac{p'}{q'}\,\geq \, K_\epsilon\, ( \frac{1}{q} +
\frac{1}{q'})^2.$$

{\bf Proof:}

Denote
$$ L= \frac{e^{-2 \, F}  }{1+ \epsilon}, \,\, M= \frac{e^{2 \, F}
}{1-\epsilon}.$$

From Corollary 1 we can suppose that:

a) first case: $L < \frac{q'}{q}<M$

or,

b) second case  $L < \frac{q}{q'}<M$.

In the first case:

$$( \frac{1}{q} +
\frac{1}{q'})^2\, < \, ( \frac{1}{q} + \frac{1}{L\, q})^2= \,
\frac{1}{q^2}\,  ( 1 + \frac{1}{L})^2$$ and,

$$ \frac{p}{q}\, - \frac{p'}{q'}\, = \frac{p \,q'\,-\,q\,p'}{q\,
q'}\,\geq \frac{1}{q\, q'}\, \geq \frac{1}{M\, q^2}.$$

Therefore,
$$\frac{p}{q}\, - \frac{p'}{q'}\,>   \frac{1}{M}\, \frac{1}{( 1 +
\frac{1}{L})^2}\,(\frac{1}{q}\, + \frac{1}{q'})^2\,.$$

The claim follows from the fact that

$$ K_\epsilon \,=\,\frac{1}{M}\, \frac{1}{( 1 +
\frac{1}{L})^2}.$$

The second case follows in a similar way.

\qed

Now we will apply the results above for the function $r(t)$.

{\bf Proposition 1:} The function $r:[a,b] \to \mathbb{R}$ is
continuous monotonous increasing. Moreover, if $   a\leq t_1 < t_2
\leq b$ and $r(t_1)\notin \mathbb{Q} $ or $r(t_2)\notin
\mathbb{Q}$, then $r(t_1) < r(t_2)$.

{\bf Proof:} For the continuity and other properties see \cite{KH}.
The rest follows from Lemma 1. \qed

 {\bf Lemma 3:} If $a\leq t_1<t_2 \leq b$, then $ g_{t_2} (x)
- g_{t_1} (x) \geq \,m\, (t_2 - t_1).$

{\bf Proof:} It follows easily from the twist condition and the
mean value theorem. \qed

{\bf Proposition 2:} There exist a set of full measure $X \subset
\mathbb{R}$ (that is, $ \mathbb{R}- X$ has Lebesgue measure zero)
and such that if $\tau \in (a,b)$ and $r(\tau) \in X$, then
$$ \limsup_{t_2 \to \tau^- , \, t_1 \to \tau^+}\, \frac{r(t_2)-r(t_1) }{  (t_2\,-\,t_1)^2
}\geq \frac{m^2}{e^{2 \, F} \, (1 + e^{2\, F} )^2}.$$

{\bf Proof:}  Given $0<\epsilon<1$, by corollary 2 there exists a
set of full measure $X_\epsilon \subset \mathbb{R}$ of irrational
numbers such that if $r(\tau) \in X_\epsilon$, then the claim of
this corollary is true for $K_\epsilon$. Let $\delta $ be such
that $( \tau - \delta , \tau + \delta) \subset (a,b)$.

The function $ g_\tau (x) - g_{\tau -\delta} (x)$ is continuous,
periodic and positive by Lemma 3. Therefore,
$$ \inf_{x \in \mathbb{R} }\, ( g_\tau (x) - g_{\tau -\delta} (x)  ) >0.$$

In the same way,
$$ \inf_{x \in \mathbb{R} }\, ( g_{\tau+ \delta}  (x) - g_{\tau} (x)  ) >0.$$

By corollary 2, there exists approximations  by excess
$\frac{p}{q}$ and by defect $\frac{p'}{q'}$ of $r(\tau)$ such that
$$\frac{1}{q'}  <    \inf_{x \in \mathbb{R} }\, ( g_\tau (x) - g_{\tau -\delta} (x)  )   , $$
$$\frac{1}{q}  <    \inf_{x \in \mathbb{R} }\, ( g_{\tau+ \delta}  (x) - g_{\tau } (x)  )   $$
and
$$  (\frac{p}{q} -\frac{p'}{q'})> K_\epsilon \, ( \frac{1}{q}\,+\,
\frac{1}{q'})^2.$$

As,
$$  \inf_{x \in \mathbb{R} }\, ( g_t (x) - g_{\tau } (x)  ) = \min_{0\leq x\leq  1 }\, ( g_t (x) - g_{\tau } (x)  )    , $$
there exist $\tau_1,\tau_2$ such that $\tau - \delta\,< \,\tau_1
\, < \tau\, < \, \tau_2\, < \, \tau + \delta$ and
$$   \inf_{x \in \mathbb{R} }\, ( g_\tau (x) - g_{\tau_1 } (x)  )
=\frac{1}{q'} , \,\mbox{and}\,\,  \inf_{x \in \mathbb{R} }\, (
g_{\tau_2} (x) - g_{\tau } (x)  ) =\frac{1}{q}.$$

Now, from Lemma 3, if  $t_1\,<\,\tau_1\, < \tau\, \, < \,
\tau_2\,<\,t_2$, then
$$   \inf_{x \in \mathbb{R} }\, ( g_\tau (x) - g_{t_1 } (x)  )
>\frac{1}{q'} , \,\mbox{and}\,\,  \inf_{x \in \mathbb{R} }\, (
g_{t_2} (x) - g_{\tau } (x)  ) >\frac{1}{q}.$$

By lemma 1,
$$ r(t_1) \leq \frac{p'}{q'}< r(\tau) < \frac{p}{q} \leq r(t_2).$$

Therefore,
$$r(t_2) - r(t_1) \geq \frac{p}{q}-\frac{p'}{q'} >  K_\epsilon \,  (\frac{1}{q}\,+\, \frac{1}{q'})^2$$

As $r$ is continuous, we get $$r(\tau_2) - r(\tau_1) \geq
K_\epsilon \, (\frac{1}{q}\,+\, \frac{1}{q'})^2.$$

Now, by Lemma 3
$$ g_{ \tau_2} (x) - g_{\tau} (x) \geq m (\tau_2-\tau), \, \, \forall x \in \mathbb{R}$$
and,
$$  \,\,g_{ \tau} (x) - g_{\tau_1} (x) \geq m\,
(\tau\,-\,\tau_1), \, \forall x \in \mathbb{R}.$$

It follows that
$$\frac{1}{q} \geq   m (\tau_2-\tau),\,\mbox{and}\,\, \frac{1}{q'} \geq   m (\tau-\tau_1).  $$

Therefore,
$$ r(\tau_2) - r(\tau_1) \geq K_\epsilon\, m^2 \, (\tau_2 -
\tau_1)^2,$$ and then

$$ \limsup_{t_1 \to \tau^- , \, t_1 \to \tau^+}\, \frac{r(t_2)-r(t_1) }{  (t_2\,-\,t_1)^2
}\geq K_\epsilon\, m^2.$$

The claim of the proposition follows from taking $$X=
\cap_{n=2}^\infty \, X_{\frac{1}{n}}.$$

\qed

\end{document}